# A modified three-point Secant method with improved rate and characteristics of convergence

Ababu T. Tiruneh[1]

[1] University of Eswatini, Department of Environmental Health Science, P.O. Box 369, Mbabane, Eswatini

**Abstract**

This paper presents a modification of Secant method for finding roots of equations that uses three points for iteration instead of just two. The development of the mathematical formula to be used in the iteration process is provided together with the proof of the rate of convergence which is 1.84 and is the same as the rate of convergence of Mueller's method of root finding. Application examples are given where it is demonstrated that for equations involving ill-conditioned cases, the proposed method has better convergence characteristics compared to Newton and Secant methods.

## 1.    Introduction

Many problems in science and engineering require solving non-linear equations or systems of equations requiring trial and error procedure because of the difficulty or impracticality of finding direct analytical solutions. For example, solutions to problems that require solving polynomial equations of degree five and higher are proven to be impossible to express in analytical terms using radicals.  Several trial and error methods and their variations have evolved over time in solving such equations. Newton's method is one such classical method that requires evaluation of the function and its derivative to estimate a linear approximation to the root. The method has convergence of order two near the root. Pathological cases may arise when the iteration approaches local extremum point where the calculation displays erratic behavior or results in singularity or divergence away from the desired root. The method also gets slow in equations that have roots of multiplicities (Gerald  and Wheatley, 1994).

Variations of Newton's method are plenty. Examples include a third order method that involves evaluation of a function and two derivatives (Weeraksoon and Fernando, 2000). Further improvement (fourth order) has been suggested by Traub requiring the same number of function and derivative evaluations (Traub, 1982). Similar fourth order convergence is claimed by Sanchez and Barrero through composite function evaluations (Sanchez and Barrero, 2011). More recently, sixth and higher order convergence have been achieved by a number of researchers (Sharma and Guha, 2007; Chun, 2007; Kou and Wang, 2007; Kou, 2007; Kou and Li, 2007; Kou, *et al*., 2009; Parhi and Gupta, 2008). The stability of Newton method in all these improvements may be an issue still relevant to explore though higher order convergence is undoubtedly a deserved merit of these methods.

The Secant method, also known as Regula falsi or the Method of cords, is another linear approximation to the root that requires two points and does not require evaluating derivatives. The order of convergence of Secant method is 1.618 and Secant method may also face converegence problems similar to Newton method. A new class of Secant like methods have recently been developed that employ more than one point of the iteration but also include evaluation of the derivative. Examples of such methods include the methods developed by Tiruneh *et al*. (Tiruneh *et al*., 2013), Fernandez-Torres (2015) and Tukral (2018[a b]).

Mueller's method is a quadratic equation approximation to the root generally involving three points of the iteration and is an extension of the Secant method as it also does not require evaluating the derivative (Mueller, 1956). Mueller's method starts with three function evaluations to begin with but continues with one additonal evaluation as the iteration progresses. The method of convergence of Mueller's ,method is of order 1.84. Mueller's method, however, can possibly converge to non-real, complex roots unless the function value of one of the points is opposite in sign to the other two ( Mekwi, 2011). Mueller's method also faces degeneration problems if the two points coincide where by the method reduces to the Secant method.

Problems of converegence of the traditional Newton and Secant methods have been tackled through approaches using a hybrid of methods. Sidi (2008) used a method that involves a multi-point Secant method whereby an n-degree polynomial is fitted using the previous points of iteration and Newton method is applied in which the first derivative of the fitted polynomial replaces the derivatiave of the actual function in the Newton formula. A method that combines bisection with that of Secant method has been suggested by Dekker (Dekker, 1969). In this method the function evaluations of bisection and Secant approaches are compared and the new point resulting in estimate of function value that is closer to the root is chosen for the next trial and error procedure. Brent (1973) suggested a procedure using root bracketing and inverse quadratic extrapolation to the root. It is an improvement over Dekker's method in terms of improving the rate of convergence. The Leap-frogging method (Kasturiarachi, 2002) uses a hybrid of Newton and Secant methods for iteration obviously resulting in improved order of converegnce to cubic convergence.

## 2. Method development

The new proposed method for root finding is an iterative technique that is based on applying Secant method to the three most recent estimates of the root. Figure 1 shows the starting point of the iteration involving the first three points. While points $(x_0, y_0)$ and $(x_1, y_1)$ may be chosen arbitrarily to start the iteration, the third starting point $(x_2, y_2)$ is better determined by application of the traditional Secant method. Such an approach of estimating the additional starting point has been proven to improve the iteration process as suggested by Thukral (1918[b]) in applying the techniques to the method developed by Tiruneh *et al* (2013).

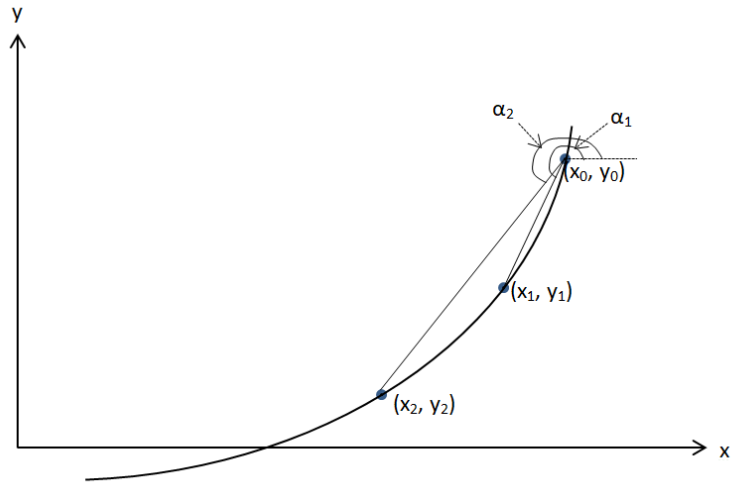

Figure 1: The starting point of the new method showing the first three points of the iteration

Referring to the x-y curve shown in Figure 1, the three distinct points are defined to be lying on the curve with coordinates $(x_0, y_0)$, $(x_1, y_1)$ and $(x_2, y_2)$. The initial point $(x_0, y_0)$ is taken as the reference point from which the angles of inclination $\alpha_i$ of lines connecting this reference point to the other two points, namely, $(x_1, y_1)$ and $(x_2, y_2)$ are measured. The variable m corresponding to the tangent of this angle of inclination for any point $(x_i, y_i)$ lying on the curve is thus defined as:

$$m_i = \tan(\alpha_i) = \frac{y_i - y_0}{x_i - x_0}$$

Analogously, the $m_1$ and $m_2$ values corresponding to the points 1 and 2 are calculated as:

$$m_1 = \tan(\alpha_1) = \frac{y_1 - y_0}{x_1 - x_0} \quad ; \quad m_2 = \tan(\alpha_2) = \frac{y_2 - y_0}{x_2 - x_0}$$

Now a new curve (m, y) is drawn with the m values replacing the corresponding x values as the independent variable. Figure 2 shows this curve.

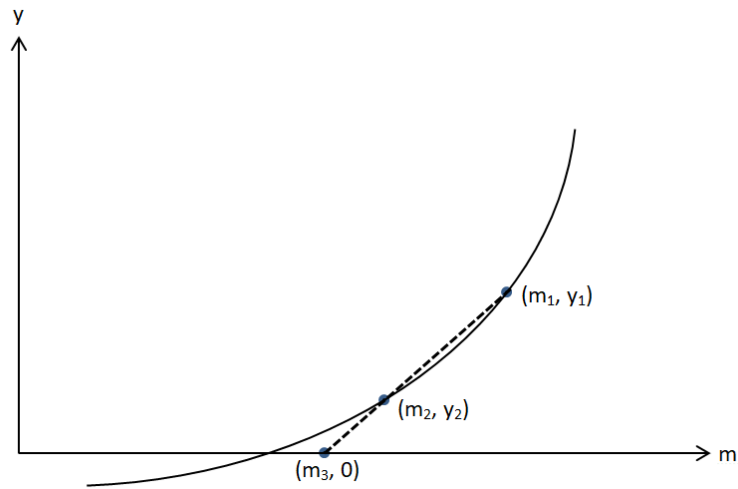

Figure 2: A new m-y curve drawn using, as the new independent variable, the tangent (m) of angle of inclination $\alpha_1$ and $\alpha_2$ of the points $(x_1, y_1)$ and $(x_2, y_2)$ measured from the reference point $(x_0, y_0)$

Now the regular Secant method is applied to this new curve shown in Figure 2 to determine the next approximation $m_3$ (and hence by extension $x_3$) to the root;

$$m_3 = m_2 - \frac{y_2}{\left(\frac{y_2 - y_1}{m_2 - m_1}\right)} \qquad (1)$$

Replacing the $m_3$ value with the corresponding $x_3$ value in the Equation (1);

$$\frac{y_3 - y_0}{x_3 - x_0} = m_2 - \frac{y_2}{\left(\frac{y_2 - y_1}{m_2 - m_1}\right)}$$

since $y_3 = 0$ as the next point of root approximation as shown in Figure 2;

$$\frac{-y_0}{x_3 - x_0} = m_2 - \frac{y_2}{\left(\frac{y_2 - y_1}{m_2 - m_1}\right)}$$

Reversing the above expression;

$$\frac{x_3 - x_0}{-y_0} = \frac{1}{m_2 - \frac{y_2}{\left(\frac{y_2 - y_1}{m_2 - m_1}\right)}}$$

$$x_3 - x_0 = \frac{-y_0}{m_2 - \frac{y_2}{\left(\frac{y_2 - y_1}{m_2 - m_1}\right)}}$$

$$x_3 = x_0 - \frac{y_0(y_2 - y_1)}{m_2(y_2 - y_1) - y_2(m_2 - m_1)}$$

Defining $\Delta y_{21} = y_2 - y_1$ and $\Delta m_{21} = m_2 - m_1$ finally gives;

$$x_3 = x_0 - \frac{y_0 \Delta y_{21}}{m_2 \Delta y_{21} - y_2 \Delta m_{21}} \qquad (2)$$

Continuing the iteration likewise, the k+1$^{th}$ estimate of the root will be;

$$x_{k+1} = x_{K-2} - \frac{y_{k-2}\Delta y_{k,k-1}}{m_k \Delta y_{k,k-1} - y_k \Delta m_{k,\ k-1}} \quad (3)$$

In terms of the x and y values of the three most recent estimates of the root, the formula in Equation (3) will eventually be written as;

$$x_{k+1} = x_{K-2} - \frac{y_{k-2}(y_k - y_{k-1})}{\left(\frac{y_k - y_{k-2}}{x_k - x_{k-2}}\right)(y_k - y_{k-1}) - y_k \left[\left(\frac{y_k - y_{k-2}}{x_k - x_{k-2}}\right) - \left(\frac{y_{k-1} - y_{k-2}}{x_{k-1} - x_{k-2}}\right)\right]} \quad (4)$$

Equation (4) will be used during the iteration process to estimate the next point of iteration from the previous most successive points of the iteration. Now the proof of order of convergence of the new method will be given.

*2.1 Proof of convergence*

Defining the error of estimate of the root r with respect to the ith estimate, $x_i$, as;

$$e_i = x_i - r$$

Equation (4) can now be written as:

$$e_{k+1} = e_{K-2} - \frac{y_{k-2}(y_k - y_{k-1})}{\left(\frac{y_k - y_{k-2}}{e_k - e_{k-2}}\right)(y_k - y_{k-1}) - y_k \left[\left(\frac{y_k - y_{k-2}}{e_k - e_{k-2}}\right) - \left(\frac{y_{k-1} - y_{k-2}}{e_{k-1} - e_{k-2}}\right)\right]} \quad .. \quad (5)$$

For points that are sufficiently close to the root $x_r$, Taylor series expansion can be used to estimate the y values in terms of the error terms;

$$y_k = (y_r = 0) + c_1 e_k + c_2(e_k)^2 + c_3(e_k)^3 + c_4(e_k)^4 + \ldots + c_n(e_k)^n + O(e_k)^{n+1} \ldots (6)$$

$$y_{k-1} = (y_r = 0) + c_1 e_{k-1} + c_2(e_{k-1})^2 + c_3(e_{k-1})^3 + \ldots + c_n(e_{k-1})^n + O(e_{k-1})^{n+1} .. (7)$$

$$y_{k-2} = (y_r = 0) + c_1 e_{k-2} + c_2(e_{k-2})^2 + c_3(e_{k-2})^3 + \ldots + c_n(e_{k-2})^n + O(e_{k-2})^{n+1} \ldots (8)$$

$$y_k - y_{k-1} = c_1(e_k - e_{k-1}) + c_2(e_k - e_{k-1})^2 + c_3(e_k - e_{k-1})^3 + c_4(e_k - e_{k-1})^4 + \ldots \ldots c_n(e_k - e_{k-1})^n \ldots \ldots (9)$$

$$\frac{y_k - y_{k-2}}{e_k - e_{k-2}} = c_1 \frac{(e_k - e_{k-2})}{e_k - e_{k-2}} + c_2 \frac{(e_k - e_{k-2})^2}{e_k - e_{k-2}} + c_3 \frac{(e_k - e_{k-2})^3}{e_k - e_{k-2}} + c_4 \frac{(e_k - e_{k-2})^4}{e_k - e_{k-2}}$$
$$+ \ldots \ldots c_n \frac{(e_k - e_{k-2})^n}{e_k - e_{k-2}} \ldots \ldots \ldots (10)$$

$$= c_1 + c_2(e_k + e_{k-2}) + c_3(e_k^2 + e_k e_{k-2} + e_{k-2}^2)$$
$$+ c_4(e_k^3 + e_k^2 e_{k-2} + e_k e_{k-2}^2 + e_{k-2}^3) + \ldots \ldots \ldots (11)$$

Similarly;

$$\frac{y_{k-1} - y_{k-2}}{e_{k-1} - e_{k-2}} = c_1 + c_2(e_{k-1} + e_{k-2}) + c_3(e_{k-1}^2 + e_{k-1}e_{k-2} + e_{k-2}^2)$$
$$+ c_4(e_{k-1}^3 + e_{k-1}^2 e_{k-2} + e_{k-1}e_{k-2}^2 + e_{k-2}^3) + \ldots \ldots \ldots (12)$$

In all of the Taylor series expansions above, the c terms in the series are given by;

$$c_n = \frac{y^n(r)}{n!}$$

Where $y_n$ is the $n^{th}$ derivative of the function y with respect to x evaluated at the root x=r.

Substitution of the expressions given in Eq. (6) to Eq. (12) in Equation (5) above and neglecting the fourth and higher order terms of the error will give after simplification;

$$e_{k+1}$$
$$= \frac{(c_2^2 - c_1 c_3) e_k e_{k-1} e_{k-2}}{c_1^2 + c_1 c_2(e_k + e_{k-1}) + c_1 c_3(e_k^2 + e_k e_{k-1} + e_{k-1}^2) + c_1 c_2 e_{k-2} + c_2^2(e_k e_{k-1} + e_k e_{k-2} + e_{k-1} e_{k-2})}$$
$$+ c_1 c_3 (e_{k-2}^2 - e_k e_{k-1})$$

The terms in the denominator containing the error terms $e_k$, $e_{k-1}$ and $e_{k-2}$ are all small enough to be neglected compared to the term dominant term $c_1^2$. This simplification gives;

$$e_{k+1} = \left(\frac{c_2{}^2 - c_1 c_3}{c_1{}^2}\right) e_k e_{k-1} e_{k-2}$$

Defining positive real terms of the sequence $S_k$ such that;

$$S_k = \frac{|e_{k+1}|}{|e_k{}^\alpha|} \ ; \ S_{k-1} = \frac{|e_k|}{|e_{k-1}{}^\alpha|} \ ; \ S_{k-2} = \frac{|e_{k-1}|}{|e_{k-2}{}^\alpha|}$$

Expressing all the error terms in terms of $e_{k-2}$;

$$|e_{k-1}| = S_{k-2} |e_{k-2}{}^\alpha|$$

$$|e_k| = S_{k-1} |e_{k-1}{}^\alpha| = S_{k-1} S_{k-2}{}^\alpha |e_{k-2}{}^{\alpha^2}|$$

$$|e_{k+1}| = S_k |e_k{}^\alpha| = S_k S_{k-1}{}^\alpha S_{k-2}{}^{\alpha^2} |e_{k-2}{}^{\alpha^3}|$$

Finally, the error ratio below converges to the constant containing the c terms that contain the derivatives of the y function evaluated at the root $x = r$;

$$\left|\frac{e_{k+1}}{e_k e_{k-1} e_{k-2}}\right| = \left|\frac{c_2{}^2 - c_1 c_3}{c_1{}^2}\right|$$

Substituting the equivalent expressions of errors containing the $e_{k-2}$ term only;

$$\left|\frac{e_{k+1}}{e_k e_{k-1} e_{k-2}}\right| = S_k S_{k-1}{}^{\alpha-1} S_{k-2}{}^{\alpha^2 - \alpha - 1} \left|e_{k-2}{}^{\alpha^3 - \alpha^2 - \alpha - 1}\right| = \left|\frac{c_2{}^2 - c_1 c_3}{c_1{}^2}\right|$$

As the iteration approaches the root, the above ratio approaches the constant on the right hand side of Equation containing the c terms. For this to be true unconditionally, the power of the error term should approach zero, i.e.,

$$\alpha^3 - \alpha^2 - \alpha - 1 = 0 \tag{13}$$

The roots of the above third degree polynomial in Equation (13) are;

$$\alpha = 1.83929 \;;\; \alpha = -0.41964 \pm 0.60629\, i$$

The order of convergence of the proposed method is therefore 1.839 which is similar to the rate of convergence of Mueller's method of root finding. However, the proposed method always converges to the real root which is not necessarily the case with Mueller's method which may converge to imaginary roots unless one of the points has its y value that is opposite in sign to the y values of the other two points. Can we say also that we have in the process discovered a family member of Mueller's method that is having linear forms?

The proposed method has, therefore, better order of convergence compared to the regular Secant method which has order of convergence of approximately 1.68. It will now be shown in the following section that the proposed method, in addition, displays better convergence characteristics for ill-conditioned cases in which either Newton or Secant method (and in some cases both) may fail to converge to the root.

## 3. Application Examples

Examples of application of the proposed three point Secant method for locating roots of several equations are discussed below. The equations tested included those that regularly converge to the roots; equations that have roots of multiplicities in which the rate of convergence is slowed down; and equations that display pathological behavior during iteration whereby application of Newton or Secant method may fail to reach convergence.

For the purpose of determining the number of iterations, a stopping criterion is used which uses the following rule:

$$|x_k - x_{k-1}| + |y_k| < 10^{-15} \qquad (14)$$

The rate of convergence of the methods to the roots is calculated with the following formula:

$$\alpha_k = \frac{Log\,(|e_{k+1}|)}{Log\,(|e_k|)} = \frac{Log\,(|X_{k+1} - r|)}{Log\,(|X_k - r|)} \qquad (15)$$

In Equation (15), $\alpha_k$ is the order of convergence of the iterative process at the $k^{th}$ iteration step, $e_k$ and $e_{k+1}$ are the errors of estimate of the root at the $k^{th}$ and $k+1^{th}$ iteration steps respectively and r is the desired root of the equation.

*3.1 Test of the proposed method for regular cases*

Table 1 shows summary of the results of the iteration processes among the three methods for equations in which the iterations regularly convergence to the root. Examination of the rate of convergence of the proposed method for regular cases show that the method generally converges at a rate of 1.84 as the iteration approaches the root value as predicted theoretically, a rate that is between Secant and Newton methods. For the equation having roots of multiplicity, namely,

$$y = (x - 2)(x + 2)^4$$

all the three methods display slow and almost linear convergence requiring greater number of iterations with Newton method relatively faster than the other two methods followed by the proposed three point Secant method. This order of rate of convergence is also similar to the order for the regular cases except for the greater number of iterations required.

It is also interesting to note that for another function that has root of multiplicity, namely,

$$y = (x - 1)^6 - 1$$

All the three methods show regular rate of convergence rate as the iteration approaches convergence near the root x= 2. However, depending on the starting point, the total number of iterations required varies differently among the different methods. The geometric mean of the rate of convergence is higher (1.64) for the proposed three point Secant method compared with Newton method (1.20) and Secant method (1.14). The total number of iterations required is also proportional to this geometric mean of convergence. The three point Secant method in this case displays superior overall rate of convergence.

In a similar pattern, for the equation:

$$y = e^{x^2 + 7x - 30} - 1$$

The geometric mean of rate of convergence for the three point Secant method is the highest (1.24) compared to Newton method (1.18) and Secant method (1.06). The total number of iterations is proportional to this geometric mean with the proposed three point Secant method requiring marginally lower number of iterations compared to Newton and Secant methods. As Table 1 shows, for this particular equation, the total number of iterations required is 26 for the three point Secant method compared to 28 for Newton method and 39 for Secant method.

Figure 3 shows comparison of the number of iterations required for finding roots of several equations among the different methods where by the proposed three point Secant method shows convergence rate between Secant and Newton methods and better overall rate of convergence for cases displaying characteristics similar to equations having roots of multiplicities as the examples above have shown.

Table 1. Comparison of result of iterations of the Three point Secant method with Newton and traditional Secant Methods.

| Function | Root | Starting points | Number of iterations required | | |
|---|---|---|---|---|---|
| | | | Secant Method | Newton Method | Three point Secant method |
| $y = x^3 + 4x^2 - 10$ | 1.365230013414100 | 0.5, 0.55, 0.6 | 10 | 8 | 9 |
| | | 0.9, 0.95, 1.0 | 10 | 6 | 8 |
| $y = [\sin(x)]^2 - x^2 + 1$ | -1.404491648215340 | -1.0, -0.975, -0.95 | 10 | 7 | 8 |
| | | -3.5, -3.25, -3.0 | 10 | 7 | 9 |
| $y = (x - 2)(x + 2)^4$ | -2.0000000000000 | -3.1, -3.05, -3.0 | 169 | 117 | 121 |
| | 2.0000000000000 | 1.4, 1.45, 1.5 | 117 | 32 | 87 |
| $y = (x - 1)^6 - 1$ | 2.00000000000000 | 1.5, 1.55, 1.6 | 25 | 17 | 9 |
| | | 2.5, 2.55, 2.6 | 12 | 8 | 8 |
| | | 3.5, 3.55, 3.6 | 15 | 11 | 11 |
| $y = \sin(x) \cdot e^x + \ln(x^2 + 1)$ | -0.603231971557215 | -0.9, -0.85, -0.8 | 10 | 7 | 7 |
| | | -0.7, -0.65, -0.6 | 8 | 5 | 5 |
| $y = e^{x^2 + 7x - 30} - 1$ | 3.000000000000000 | 4.0, 4.05, 4.1 | 29 | 20 | 20 |
| | | 4.4, 4.45, 4.5 | 39 | 28 | 26 |
| $y = x - 3\ln(x)$ | 1.857183860207840 | 2.0, 2.05, 2.1 | 9 | 5 | 7 |
| | | 0.4, 0.45, 0.5 | 11 | 8 | 8 |

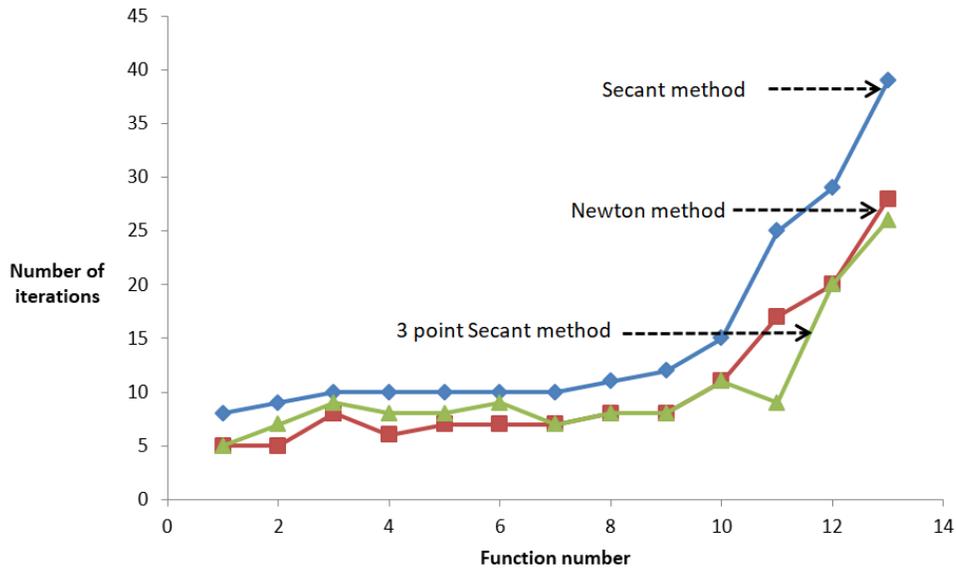

Figure 3: Comparison of number of iterations required for convergence among the different methods

## *3.2 Test of the proposed method for divergent cases*

The proposed new method has been examined for its convergence for pathological cases whereby the iteration process displays divergence or failure when applying either the Newton or Secant methods. Table 2 summarizes the results of the iteration for these pathological cases in which comparison is method between the proposed three point Secant method and the traditional Newton and Secant methods. The Table shows that, whereas either Newton or Secant method (and in some cases both) fail to converge to the root; the proposed method almost always converges to the root. This result shows that the proposed three point Secant methods has a relative advantage for ill-conditioned cases where application of either Newton or Secant or both methods may not lead to convergence.

Figures 4-6 show the characteristics of convergence for the three different methods when applied to the equation:

$$y = x^{1/3}$$

Where the root x=0 is a straightforward solution. It is seen in Figure 4 that applying Newton method results in oscillatory divergence to infinity. Similarly, application of Secant method to the same equation leads to oscillation between finite points in which the iteration fails to converge to the root as shown in Figure 5. By contrast the proposed three point Secant method shows oscillatory convergence towards the root as Figure 6 shows.

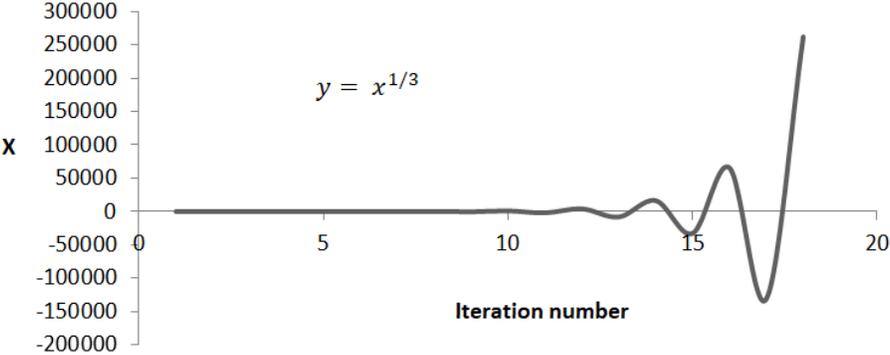

Figure 4: Divergence of iteration by Newton method for the function y = x$^{1/3}$

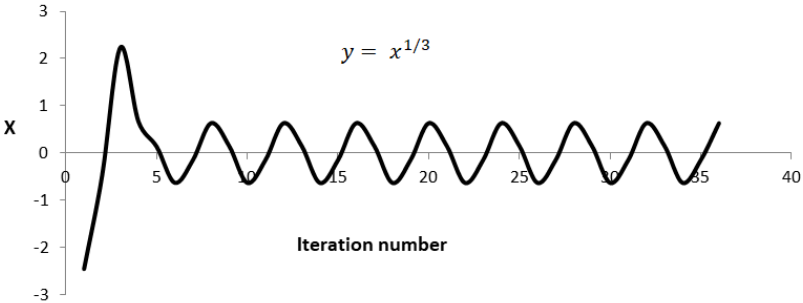

Figure 5: Oscillation of the iteration by Secant method for the function y = x$^{1/3}$

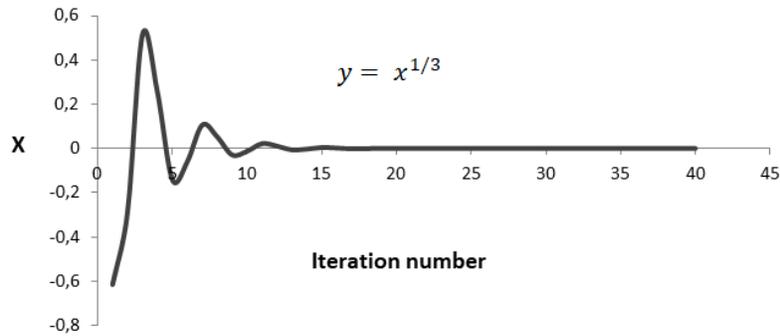

Figure 6: Oscillating convergence of the iteration by three point secant method for y= x$^{1/3}$

Another example of the advantage of the proposed method for ill-conditioned cases is the five degree polynomial equation:

$$y = x^5 - x + 1$$

Application of Newton method results in oscillation where by the iteration fails to converge to the root as shown graphically in Figure 7. In a similar pattern, application of the regular Secant method results in a gradually diverging oscillation where the method fails to converge to the root as shown in Figure 8. By contrast, application of the proposed three point Secant methods shows convergence with a rate of 1.84 as proven theoretically near the root and with an overall rate of convergence (geometric mean) of 1.15 for the starting points given in Table 2. The process of convergence for this polynomial equation as the iteration proceeds is shown graphically in Figure 9.

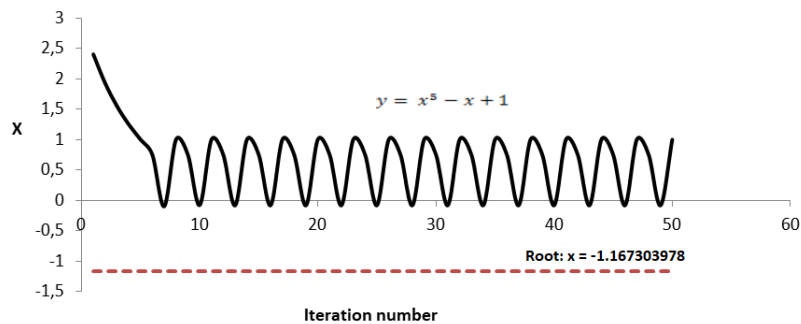

Figure 7: Constant oscillation of the iteration by Newton method for y= x$^5$- x +1

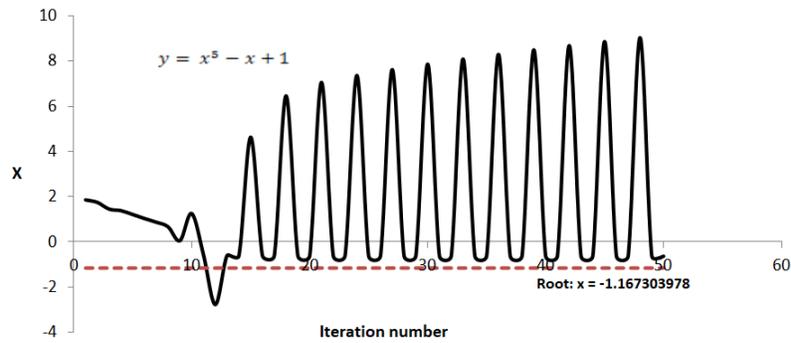

Figure 8: Oscillating divergence of the iteration by Secant method for y= $x^5$- x +1

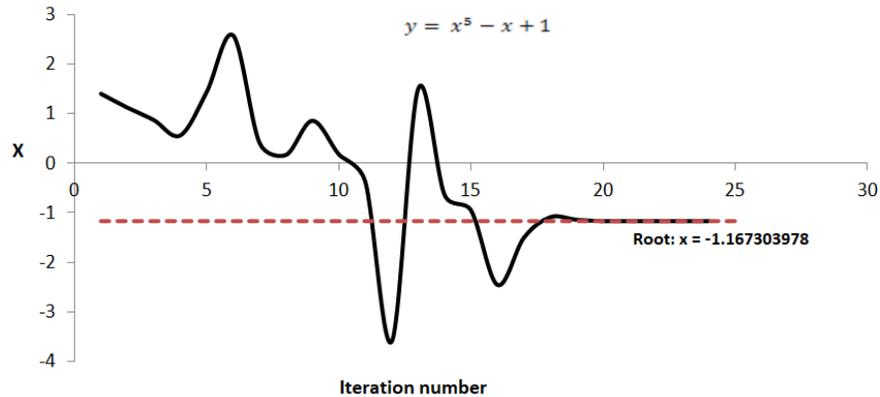

Figure 9: Oscillating convergence of the iteration by three point Secant method for y= $x^5$- x +1

## 4. Conclusion

A mathematical formula involving iteration for finding roots of non-linear equations has been developed and presented in this paper. The proposed formula is a modification of the Secant method employing three successive points of the iteration instead of just two. The method applies the Secant method to the tangent of the relative angle of inclination of the lines joining the furthest point with the two most recent points of iteration. The rate of convergence of the proposed method is of the order 1.83929 which is equivalent to the rate of convergence of Mueller's method of root finding. In fact it can be said that this method is a discovery of a variant or a family of Mueller's method that is linear in form as it is recalled that Mueller's method is also an iterative procedure that is based on three successive points of the iteration.

Table 2. Results of application of the proposed three point Secant method for cases in which either Newton or Secant method (or both) fail to converge.

| Function | Root | Comparison of number of iterations required | | | |
|---|---|---|---|---|---|
| | | Starting points | Secant Method | Newton Method | Three point Secant method |
| $y = -x^4 + 3x^2 + 2$ | 1.887207676120680 | 1, 1.5, 2.0 | 11 | Oscillates | 10 |
| | | 0.5, 0.55, 0.6 | 24 | Oscillates | 23 |
| $y = \log(x)$ | 1.000000000000000 | 3.0, 3.25, 3.5 | Fails | Fails | 9 |
| $y = Arctan(x)$ | 0.0000000000000 | 3.0, 3.25, 3.5 | Diverges | Diverges | 9 |
| | | -3.0, -3.25, -3.5 | Diverges | Diverges | 10 |
| $y = x^5 - x + 1$ | -1.167303978261420 | 2.0, 2.5, 3.0 | Oscillates | Oscillates | 27 |
| | | 7.0, 7.5, 8.0 | 112 | Oscillates | 44 |
| $y = 0.5x^3 - 6x^2 + 21.5x - 22$ | 1.7639320225002100 | 2.0, 2.5, 3.0 | 10 | Oscillates | 9 |
| | 6.236067977499790 | 5.0, 5.5, 6.0 | 10 | Oscillates | 11 |
| $y = x^{1/3}$ | 0.00000000000000 | 1.0, 1.25, 1.5 | Oscillates | Diverges | 95 |
| | | -1.0, -1.25, -1.5 | Oscillates | Diverges | 93 |
| $y = 10xe^{-x^2} - 1$ | 1.679630610428450 | 3.0, 3.25, 3.5 | Diverges | Diverges | 14 |
| | 0.101025848315685 | -1.0, -1.5, -2.0 | Diverges | Diverges | 19 |

However, the proposed method is linear in nature being based on Secant method of root finding and as such does not lead to imaginary roots unlike Mueller's method which is based on quadratic solution to the root approximation. Stated in another way, the method does not require that the function values of one of the iteration points be opposite in sign to the others in order to avoid imaginary root values during the iteration.

Examples of application of the proposed method of root finding to a variety of equations has been presented and compared with the results of iterations of Newton and Secant methods. For equations leading to regular convergence, it is shown that the proposed method has rate of convergence that lies between Secant and Newton methods as is also supported by the mathematical proof of the rate of convergence presented in this paper. For equations that display roots of multiplicities or characteristics similar to roots of multiplicities, the proposed method generally displays better overall rate of convergence compared to Secant and Newton methods.

For ill-conditioned cases in which Newton and Secant methods may fail to converge displaying oscillation, divergence to infinity or off-shooting to undesirable or invalid domain, the proposed three point Secant method almost always leads to convergence as is demonstrated by a number of examples presented in this paper. This is the inherent advantage of the proposed method over the traditional Newton and Secant methods that display pathological behavior for ill-conditioned cases. The application examples demonstrate that the proposed method has better convergence characteristics for such ill-conditioned cases compared to Newton and Secant methods.

**References**

1   Gerald, C.F. and Wheatley, P.O. (1994) *Applied numerical analysis*. Fifth Edition, 40-75.

2   Weerakoon, S. and Fernando, T.G.I. (2000) A variant of Newton's method with accelerated third order convergence, *Applied Mathematics Letters*, **13**, 87-93.

3   Traub, J.F. (1982) *Iterative methods for the solution of equations*, Prentice-Hall, Englewood.

4   Sanchez, M.G. and Barrero, J.L.D. (2011) A Technique to composite a modified Newton's method for solving nonlinear Equations, arXiv: 1106.0996v1, Cornell University Library.

5   Sharma, J.R. and Guha, R.K. (2007) A family of modified Ostrowski methods with accelerated sixth order convergence, *Appl. Math. Comput*. **190**, 111-115

6   Chun, C. (2007) Some improvements of Jarratt's method with sixth-order convergence, *Appl. Math. Comput*. **190**, 1432- 1437.

7   Kou, J. and Wang, X (2007) Sixth-order variants of Chebyshev- Halley methods for solving non-linear equations, *Appl. Math. Comput.* , **190**, 1839-1843.

8   Kou, J. (2007) On Chebyshev-Halley methods with sixth-order convergence for solving non-linear equations, *Appl. Math. Comput.*, **190**, 126-131.

9   Kou J., and Li, Y. (2007) An improvement of Jarratt method, *Appl. Math. Comput*. **189**, 1816-1821.

10  Kou, J. Li, Y.and Wang, X. (2009) Some modifications of Newton's method with fifth order convergence, *Journal of Computational and Applied Mathematics,* **209**, 146-152.

11  Parhi, S.K. and Gupta, D.K. (2008) A sixth order method for nonlinear equations, *Applied Mathematics and Computation* **203**, 50-55.

12  David E. Muller, (1956) A Method for Solving Algebraic Equations Using an Automatic Computer, Mathematical Tables and Other Aids to Computation, 10, 208-215.

13  Mekwi, W.R. (2001) Iterative methods for roots of polynomials, University of Oxford. M.Sc. thesis,

14  Dekker, T.J. (1969) Finding a zero by means of successive linear interpolation, In B. Dejon and P.Henrici (eds), Constructive Aspects of the Fundamental Theorem of Algebra, Wiley-Interscience, London, SBN 471-28300-9.

15  Brent, R.P. (1973) Algorithms for Minimization without Derivatives, Chapter 4. Prentice-Hall, Englewood Cliffs, NJ. ISBN 0-13-022335-2.

16  Kasturiarachi, A. B. (2002) A Leapfrogging Newton's method, *International Journal of Mathematical Education in Science and Technology*, **33**(4), 521-527.

17  Sidi, A. (2008). Generalization of the secant method for nonlinear equation. *Applied Mathematics E-Notes*, 8, 115-123.

18  Tiruneh, A.T., Ndlela W. N. and Nkambule S. J. (2013) A two-point Newton method suitable for non-convergent cases and with super-quadratic convergence, *Adv. Numer. Anal.* art. ID687382.

19  Fernandez-Torres, G.A (2015). A novel geometric modification of the Newton-secant method to achieve convergence of order $1+\sqrt{2}$ and its dynamics, *Mod. Sim. Eng.*, art. ID 502854.

20  Thukral, R. (2018[a]). A New Secant-type method for solving nonlinear equations. *American Journal of Computational and Applied Mathematics*, **8**(2): 32-36.

21  Thukral, R. (2018[b]) Further development of secant-type methods for solving nonlinear equations. *International Journal of Advances in Mathematics, 38(5), 45-53.*